\documentclass[11pt,reqno]{amsart}
\usepackage{yhmath}
\usepackage{fancyhdr}
\usepackage{times}
\usepackage[T1]{fontenc}
\usepackage{mathrsfs}
\usepackage{latexsym}
\usepackage[dvips]{graphics}
\usepackage{epsfig}
\usepackage{hyperref, flafter}
\usepackage{amsmath,amsfonts,amsthm,amssymb,amscd}
\usepackage{color}
\usepackage{geometry}               
\usepackage{epstopdf}
\usepackage{verbatim}

\allowdisplaybreaks

\newtheorem{thm}{Theorem}[section]

\newtheorem{lem}[thm]{Lemma}
\newtheorem*{thm*}{Theorem}
\newtheorem*{con*}{Conjecture}
\newtheorem*{lem*}{Lemma}
\newtheorem{prop}[thm]{Proposition}

\def\g{\gamma}
\def\r{\rho}
\def\L{\log T}

\begin{document}
	
	\title[log derivative of zeta]
	{Mean values of the logarithmic derivative of the Riemann zeta-function near the critical line}
	
	\author{Fan Ge}

	\address{Department of Mathematics, William \& Mary, Williamsburg, VA, United States}
	
	\email{ge@wm.edu}
	
	\baselineskip=15pt
	
	\begin{abstract}
		Assuming the Riemann Hypothesis and a hypothesis on small gaps between zeta zeros (see equation~\eqref{eq ES 2k} below for a precise definition), we prove a conjecture of Bailey, Bettin, Blower, Conrey, Prokhorov, Rubinstein and Snaith~\cite{bailey2019mixedmoments}, 
		which states that for any positive integer $K$ and real number $a>0$,
		\begin{align*}
			\lim_{a \to 0^+}\lim_{T \to \infty}  \frac{(2a)^{2K-1}}{T (\log T)^{2K}}	\int_{T}^{2T} \left|\frac{\zeta'}{\zeta}\left(\frac{1}{2}+\frac{a}{\log T}+it\right)\right|^{2K} dt =  \binom{2K-2}{K-1}.
		\end{align*}
		When $K=1$, this was essentially a result of Goldston, Gonek and Montgomery~\cite{goldston2001pair} (see equation~\eqref{eq ggm} below).
	\end{abstract}
	
	\maketitle
	
	\section{Introduction}
	
	Throughout we assume the Riemann Hypothesis (RH). 
	Let
	\begin{align*}
		\mathcal I_K(a, T)=	\int_{T}^{2T} \left|\frac{\zeta'}{\zeta}\left(\frac{1}{2}+\frac{a}{\log T}+it\right)\right|^{2K} dt .
	\end{align*}
	In an important paper of Goldston, Gonek and Montgomery~\cite{goldston2001pair}, they studied $\mathcal I_1(a, T)$
	and find its close connections to Montgomery's Pair Correlation Conjecture (see Montgomery~\cite{montgomery1973pair}) and to the distribution of prime numbers in short intervals (see also Selberg~\cite{selberg1943prime}) when $a$ is of constant size. Generalizations of $\mathcal I_K(a, T)$ for larger $K$ have also been studied by Farmer, Gonek, Lee and Lester in~\cite{farmer2013logderiv}, in which they related these integrals to higher correlation functions of zeta zeros.
	
	In~\cite{goldston2001pair} Goldston et al. also proved that 
	\begin{align}\label{eq ggm}
		\mathcal I_1(a,T) \sim \frac{1}{2a} T(\L)^2 \quad \text{ as } a =a(T) \to 0,
	\end{align}
	assuming RH and the Essential Simplicity hypothesis (ES). To state (ES), let $\rho=1/2 +i\g$ be a generic zero of zeta, and let
	\begin{align*}
		N^*(T) = \sum_{0<\g<T} m_\g
	\end{align*}
	where $m_\g$ is the multiplicity of $\r$ and the sum counts multiplicity of zeros. We also define
	\begin{align*}
		N(\beta, T) =  \sum_{\substack{0<\g, \g'<T \\ 0< \g-\g' < 2\pi\beta/\L}} 1.
	\end{align*}
	Recall that 
	\begin{align*}
		N(T)=\sum_{0<\g<T} 1 \sim  \frac{T}{2\pi} \L.
	\end{align*}
	Then, the (ES) states that
	\begin{align}\label{eq ES}\tag{ES}
		N^*(T) \sim \frac{T}{2\pi} \L, \quad \text{ and } \quad  N(\beta, T) = o (T\L)  \text{ for any } \beta=\beta(T) \to 0.
	\end{align}
	
	It is widely believed that the circular unitary ensemble (CUE) models the Riemann zeta-function regarding zero statistics, value distribution, and more (for example, see Montgomery~\cite{montgomery1973pair}, Keating and Snaith~\cite{keating2000momentszeta}). Bailey, Bettin, Blower, Conrey, Prokhorov, Rubinstein and Snaith~\cite{bailey2019mixedmoments} proved asymptotics for the analogue of $\mathcal I_K(a, T)$ when $a \to 0$ in the setting of CUE for all positive integers $K$. The expressions in these asymptotics lead them to conjecture that
	\begin{align*}
		\lim_{a \to 0^+}\lim_{T \to \infty}  \frac{(2a)^{2K-1}}{T (\log T)^{2K}} \cdot	\mathcal I_K(a, T) =  \binom{2K-2}{K-1}.
	\end{align*}
	When $K=1$, this conjecture agrees with the result of Goldston et al~\eqref{eq ggm} apart from a small difference in how the limits are taken. 
	
	The purpose of this paper is to prove the above conjecture of Bailey et al., assuming RH and a hypothesis which we call the $2K$-tuple Essential Simplicity hypothesis (ES $2K$).  To state (ES $2K$), let  $\r_1, \r_2, ... $ be the list of zeta zeros on the upper half plane ordered with heights and counted with multiplicity, and define
	%
	%	$N(K,v, T)  $ to be the number of $2K$-tuples $(\lambda_1, ..., \lambda_{2K}) $ of zeta zeros with height in $[T, 2T]$ such that all of $\lambda_1, ..., \lambda_{2K}$ lie in an interval of length $v/\L$, and such that not all of  $\lambda_1, ..., \lambda_{2K}$ have the same label $\r_i$ in the list $\r_1, \r_2, ...$ of zeta zeros.
	%
	$N(K,v, T)  $ to be the number of $2K$-tuples $(\lambda_1, ..., \lambda_{2K}) $ of zeta zeros satisfying
	\begin{itemize}
		\item $\lambda_j \in [T, 2T]$ for all $j=1,..., 2K$;
		\item $|\lambda_i-\lambda_j| \le \frac{v}{\L}$ for all $i, j=1,..., 2K$;
		\item not all of  $\lambda_1, ..., \lambda_{2K}$ have the same label $\r_i$ in the list $\r_1, \r_2, ...$ of zeta zeros.
	\end{itemize}
	(If all the zeros of zeta are simple, then the last condition is the same as saying that not all of $\lambda_1, ..., \lambda_{2K}$ are the same complex number.)
	Then, our (ES $2K$) states that
	\begin{align}\label{eq ES 2k}\tag{ES $2K$}
		\lim_{T \to \infty}   \frac{N(K, v, T)}{N(T)} \to 0 \quad \text{ as } v \to 0.
	\end{align}
	Therefore, roughly speaking, this hypothesis states that there are not too many small gaps between zeta zeros. For $K=1$, this assumption is essentially the same as the (ES) assumed in the result of Goldston et al.~\eqref{eq ggm}, apart from a small difference in how we take the limits. (If we let $v=v(T)$ depend on $T$ in our assumption as in (ES), our result below modifies accordingly to a generalization of~\eqref{eq ggm}.) Note also that our~\eqref{eq ES 2k} is consistent with, but much weaker than, the more general GUE hypothesis for $n$-correlations of zeta zeros. 
	
	\begin{thm}\label{thm}
		Let $K$ be a positive integer, and $a\in \mathbb R$. Assume RH and (ES $2K$).  We have
		\begin{align*}
			\lim_{a \to 0^+}\lim_{T \to \infty}  \frac{(2a)^{2K-1}}{T (\log T)^{2K}}	\int_{T}^{2T} \left|\frac{\zeta'}{\zeta}\left(\frac{1}{2}+\frac{a}{\log T}+it\right)\right|^{2K} dt =  \binom{2K-2}{K-1}.
		\end{align*}
	\end{thm}
	
	We remark that the method used in proving~\eqref{eq ggm} (i.e., the case when $K=1$) in Goldston et al.~\cite{goldston2001pair} does not seem to generalize to $K\ge 2$.
	The method in the work of Farmer et al.~\cite{farmer2013logderiv} could possibly be used to obtain estimates for $\mathcal I_K(a, T)$ (and for more general integrals), assuming certain knowledge on the correlation functions of zeta zeros; however, according to~\cite{farmer2013logderiv}, the calculation seems quite difficult in general. Thus, it is currently not clear whether we could  obtain Theorem~\ref{thm} using this method. Another possible approach to studying $\mathcal I_K(a, T)$ is via the ratios conjectures for zeta (see Conrey, Farmer and Zirnbauer~\cite{conrey2008autocorrelationofratios}, Conrey and Snaith~\cite{conrey2008correlations}). Assuming the ratios conjectures, it is possible to obtain  estimates (also involving complicated calculation) for $\mathcal I_K(a, T)$  for each $K$ individually; however, it seems that obtaining a formula for a general $K$ is currently unknown;   see~\cite{bailey2019mixedmoments} for relevant discussion in the CUE case.  
	On the contrary, our approach here is to directly calculate the leading asymptotics of the integral when $a\to 0$, instead of working on the asymptotics of $\mathcal I_K(a, T)$ first for $a$ of constant size (which is more important but more difficult) by appealing to its relations to the correlation functions or ratios conjectures. Moreover, our proof reveals transparently where the binomial coefficient in the result comes from. 
	
	We also remark that our result shows  the limiting asymptotic behavior of $\mathcal I_K(a, T)$ when $a\to 0$ (as expected in Theorem~\ref{thm}) only requires that there are not too many very small gaps between zeta zeros, but \textit{does not require  finer information on the distribution of typical sizes of zero gaps}. Thus, in particular, we expect that the leading term of $\mathcal I_K(a, T)$ as $a\to 0$ cannot distinguish the believed distribution of zeta zeros from the Alternative Hypothesis. (See Baluyot~\cite{baluyot2016ah} for a discussion of $\mathcal I_1(a, T)$ under the Alternative Hypothesis; see also Ki~\cite{ki2008zeros} for an estimate of $\mathcal I_1(a, T)$ under an assumption on zeros of $\zeta'(s)$, as well as Soundararajan~\cite{soundararajan1998horizontal}, Zhang~\cite{zhang2001zeros}, and Ge~\cite{ge2017gaps} for implication of Ki's assumption to small gaps between zeta zeros.)
	
	In the next section we state some propositions  and deduce the theorem from them. We prove these propositions in  subsequent sections.
	
	\section{Propositions and Proof of Theorem~\ref{thm}}
	
	Our first proposition, which is Lemma 5 in Radziwi{\l}{\l}~\cite{radziwill2014gaps}, plays a key role in our proof. 
	\begin{prop} \label{lem rad}
		Assume RH. Let $0<c<1$. Uniformly in $t\in [T, 2T]$, $N\le T$, and $\frac{1}{2}\le \sigma \le \frac{1}{2}+\frac{1}{\log N}$, we have
		\begin{align}\label{eq lem rad}
			\frac{\zeta'}{\zeta}(s) = \sum_{|s-\r|<\frac{c}{\L}} \frac{1}{s-\r} +O\left(\frac{\L}{c}\cdot \mathcal E_{T, N}(s)\right)
		\end{align}
		where $s=\sigma + it$,
		\begin{align}\label{eq E}
			\mathcal E_{T, N}(s) = \frac{\L}{\log N} + \frac{1}{\log N} \cdot \left( \left|A_N\left(\frac{1}{2}+\frac{2}{\log N}+it\right)\right| + \left| B_N\left(\frac{1}{2} + it \right)\right|\right),
		\end{align}
		with 
		\begin{align*}
			A_N(s)=\sum_{n\le N} \frac{\Lambda (n) W_N(n)}{n^s}, \quad  \quad 	B_N(s)=\sum_{n\le N} \frac{\Lambda (n)}{n^s}\left(1-\frac{\log n}{\log N}\right) ,
		\end{align*}
		and $W_N(n)=1$ for $1\le n \le \sqrt{N}$ and  $W_N(n)=1-\frac{\log n}{\log N}$ for $ \sqrt{N} < n \le N$.
	\end{prop}
	
	Thus, roughly speaking, for $s$ close to the critical line the value of $\zeta'/\zeta(s)$ can be approximated using zeros very nearby $s$ with a typically controllable error.  Such type of results was originally proved by Selberg~\cite{selbergzeta'} for $c=1$ and $s$ on the critical line. (Selberg actually proved an unconditional version.) The key ingredient in the proof is a bound of $\zeta'/\zeta(s_0)$ for a suitable $s_0$ to the right of the critical line, and this was obtained by Selberg in an ingenius way. See also  Ge~\cite{ge2022logarithmic} for a CUE analogue of Proposition~\ref{lem rad}.

	From this point on,  we let 
	\begin{align}\label{eq def s}
		s=\frac{1}{2}+\frac{a}{\log T}+it
	\end{align} 
	with $t\in [T, 2T]$ and $N=\big\lfloor \big(\frac{T}{\L}\big)^{1/K} \big \rfloor$. We also let $c$ be such that \begin{align}\label{eq def c}
		|s-\r|<\frac{c}{\L} \iff |t-\g| < \frac{b}{\L}
	\end{align} where $b=a^\delta$ with $\delta=1/4$, say. Thus, both $b$ and $c$ depend on $a$. As we eventually take the limit $a\to 0$, there is no harm to think that $c$ is about the same size as $b$.
	
	\begin{prop} \label{lem E int}
		With such $s$ and $N$, and with $\mathcal E_{T,N}(s)$ defined in \eqref{eq E}, we have
		\begin{align*}
			\int_{T}^{2T} \left| \mathcal E_{T, N}(s)\right|^{2K} dt  \ll_K T.
		\end{align*}
	\end{prop}

	\begin{prop} \label{lem zero int}
		Assume RH and (ES $2K$). With $s$ in~\eqref{eq def s} and $c$ in~\eqref{eq def c}, we have
		\begin{align*}
			\lim_{a \to 0^+} \lim_{T \to \infty}  \frac{(2a)^{2K-1}}{T (\log T)^{2K}}	\int_{T}^{2T} \left| \sum_{|s-\r|<\frac{c}{\L}} \frac{1}{s-\r} \right|^{2K} dt  =  \binom{2K-2}{K-1}.
		\end{align*}
	\end{prop}
	
	\textit{Proof of Theorem~\ref{thm}.} 
	We use Proposition~\ref{lem rad} with $c$ in~\eqref{eq def c} to expand the $2K$-th power of $\zeta'/\zeta$ into a sum of products, and then use H\"older's inequality for each summand. With
	Proposition~\ref{lem zero int} we see that it suffices to prove
	\begin{align*}
		\lim_{a \to 0^+} \lim_{T \to \infty}  \frac{(2a)^{2K-1}}{T (\log T)^{2K}}	\int_{T}^{2T} \left| E \right|^{2K} dt = 0,
	\end{align*}
	where $E$ is the error term in \eqref{eq lem rad}. From Proposition~\ref{lem E int} it follows immediately that
	\begin{align*}
		\int_{T}^{2T} \left| E \right|^{2K} dt  \ll_K \frac{T(\L)^{2K}}{c^{2K}}.
	\end{align*}
	Since $c \gg a^{1/4}$, we see that $\frac{1}{c^{2K}} = o\big(\frac{1}{a^{2K-1}}\big)$ as $a\to 0^+$. The theorem is proved. \qed

	\section{Proof of Proposition~\ref{lem E int}}

	The treatment in this section is similar to some moments computations in Radziwi{\l}{\l}~\cite{radziwill2014gaps}. We  require the following Majorant Principle, which is Theorem 3 in Chapter 7 of Montgomery's book~\cite{montgomery1994tenlectures}.
	
	\begin{lem}\label{lem majorant} 
		Let $x_1, ..., x_N$ be real numbers, and suppose that $|c_n| \le C_n$ for all $n$. Then
		\begin{align*}
			\int_{-T}^T \left|\sum_{n=1}^{N} c_n e(x_n t)\right|^2 dt \le 3 	\int_{-T}^T \left|\sum_{n=1}^{N} C_n e(x_n t)\right|^2 dt .
		\end{align*}
		Here, as usual, $e(x)=e^{2\pi i x}$.
	\end{lem}
	
	We also need a familiar result of Soundararajan (Lemma 3 in~\cite{soundararajan2009moments}) for moments of Dirichlet polynomials.
	
	\begin{lem}\label{lem sound} 
		Let $T$ be large, and let $2 \le N \le T$. Let $k$ be a natural number such that $N^k  \le T / \log T$. For any complex numbers $a(p)$ where $p$ is a prime number, we have
		\begin{align*}
			\int_{T}^{2T} \left|\sum_{p\le N} \frac{a(p)}{p^{1/2+it}}\right|^{2k} dt \ll k!T\left(\sum_{p\le N} \frac{|a(p)|^2}{p} \right)^k.
		\end{align*}
	\end{lem}
	
	Now we are ready to prove Proposition~\ref{lem E int}.
	
	Expanding the $2K$-th power and using H\"older's inequality repeatedly, we see that it suffices to prove the following bounds: 
	\begin{align*}
		&	\int_{T}^{2T} \left| \frac{\L}{\log N} \right|^{2K} dt  \ll_K T,  \quad	\int_{T}^{2T} \left| \frac{1}{\log N} \cdot A_N\left(\frac{1}{2}+\frac{2}{\log N}+it\right) \right|^{2K} dt  \ll_K T, \\
		& \qquad \qquad \text{  and  }  	\int_{T}^{2T} \left| \frac{1}{\log N} \cdot B_N\left(\frac{1}{2} + it \right) \right|^{2K} dt  \ll_K T.
	\end{align*}
	
	The first inequality is clear, since $\frac{\L}{\log N} \ll_K 1$. For the second inequality, we first apply the Majorant Principle to the Dirichlet polynomial $A_N^K$ and obtain that
	\begin{align*}
		\int_{T}^{2T} \left| A_N\left(\frac{1}{2}+\frac{2}{\log N}+it\right) \right|^{2K} dt &  \ll 	\int_{-2T}^{2T} \left| A_N\left(\frac{1}{2}+\frac{2}{\log N}+it\right) \right|^{2K} dt \\
		& \ll 	\int_{-2T}^{2T} \left| \sum_{n=1}^{N} \frac{\Lambda(n)}{n^{1/2+it}} \right|^{2K} dt .
	\end{align*}
	This last integral is equal to
	\begin{align*}
		\int_{-2T}^{2T} \left| \sum_{p\le N} \frac{\log p}{p^{1/2+it}} + \sum_{p^2\le N} \frac{\log p}{p^{1+i2t}} + O(1) \right|^{2K} dt ,
	\end{align*}
	where we can again expand the power, use H\"older's inequality, and then apply Lemma~\ref{lem sound} to bound the integral of the $2K$-th power of each of the two sums over $p$ in above. Using standard estimates, we obtain an upper bound which is $\ll_K T(\log N)^{2K}$.

	The third inequality is proved in a similar way.
	\qed
	
	\section{Proof of Proposition~\ref{lem zero int}}
	
	Recall that $s, b,$ and $c$ are defined by~\eqref{eq def s} and ~\eqref{eq def c}. We separate the interval $[T, 2T]$ into three subsets, as follows. Let
	\begin{align*}
		&	\mathcal T_0 = \left\{ t\in [T, 2T]: \text{ there is no zero } \r \text{ with } |s-\r| < \frac{c}{\L} \right\} \\
		&  S_1 = \left\{ \r: \g\in [T, 2T], \text{ there is no zero } \r^*\ne \r \text{ with } |\r^*-\r| < \frac{2b}{\L} \right\} \\
		&  \mathcal T_1 =\left \{ t\in [T, 2T]: \text{ there is a zero } \r \in S_1 \text{ with } |s-\r| < \frac{c}{\L} \right\} \\
		& \text{ and } \mathcal T_2 = [T, 2T] - \mathcal T_0 - \mathcal T_1.
	\end{align*}
	Therefore, we have
	\begin{align*}
		&	\lim_{a \to 0^+} \lim_{T \to \infty}  \frac{(2a)^{2K-1}}{T (\log T)^{2K}}	\int_{T}^{2T} \left| \sum_{|s-\r|<\frac{c}{\L}} \frac{1}{s-\r} \right|^{2K} dt  \\
		=  & \lim_{a \to 0^+} \lim_{T \to \infty}  \frac{(2a)^{2K-1}}{T (\log T)^{2K}}	\left(	\int_{\mathcal T_0} + 	\int_{\mathcal T_1} + 	\int_{\mathcal T_2}\right) \left| \sum_{|s-\r|<\frac{c}{\L}} \frac{1}{s-\r} \right|^{2K} dt  
	\end{align*}
	We estimate the integrals over $\mathcal T_i$ for $i=0,1,2$. Trivially we have $\int_{\mathcal T_0}=0$ since the sum in the integrand is empty by the definition of $\mathcal T_0$.
	
	Note that according to~\eqref{eq def c}, for $t\in  \mathcal T_1$ there is exactly one zero $\r$ with $ |s-\r| < \frac{c}{\L} $. We denote such zero by $\r_s$. Observe that we can also write 
	\begin{align*}
		\mathcal T_1 = [T, 2T] \ \bigcap \left( \bigcup_{\r\in S_1} \left\{s:  |s-\r| < \frac{c}{\L}\right\} \right).
	\end{align*}
	Thus, 
	\begin{align*}
		\int_{\mathcal T_1} \left| \sum_{|s-\r|<\frac{c}{\L}} \frac{1}{s-\r} \right|^{2K} dt  & = 	\int_{\mathcal T_1} \frac{1}{|s-\r_s|^{2K}} dt \\ 
		& = (1+o(1)) \cdot \sum_{\r\in S_1} \int_{|s-\r|<\frac{c}{\L}} \frac{1}{|s-\r|^{2K}} dt 
	\end{align*}
	where the $o(1)$ accounts for the negligible errors  from zeros near the boundary (i.e., near $T$ or $2T$). For each integral above we extend the range of integration to $\mathbb R$ with small error, as follows. We write $s=\sigma + it$ and
	\begin{align*}
		\int_{|s-\r|<\frac{c}{\L}} \frac{1}{|s-\r|^{2K}} dt & = \int_{|s-\r|<\frac{c}{\L}} \frac{1}{\left((\sigma - 1/2)^2 + (t-\g)^2 \right)^{K}} dt \\
		& = \int_{-\infty}^{\infty} \frac{1}{\left((\sigma - 1/2)^2 + (t-\g)^2 \right)^{K}} dt  + O\left(\int_{\g+\frac{b}{\L}} \frac{1}{\left((\sigma - 1/2)^2 + (t-\g)^2 \right)^{K}} dt \right)\\
		& =: I + O(J),
	\end{align*}
	say. It is straightforward to compute (or using Mathematica) that 
	\begin{align*}
		I = \frac{1}{(\sigma-1/2)^{2K-1}}\cdot \sqrt{\pi} \cdot \frac{\Gamma(K-1/2)}{\Gamma (K)}
	\end{align*}
	and that
	\begin{align*}
		J = \frac{1}{(\sigma-1/2)^{2K-1}}\cdot \int_{1/b}^\infty \frac{1}{(1+t^2)^K} dt = o(I)
	\end{align*}
	as $a\to 0$.
	Therefore, 
	\begin{align*}
		\int_{\mathcal T_1} \left| \sum_{|s-\r|<\frac{c}{\L}} \frac{1}{s-\r} \right|^{2K} dt  & = 	(1+o(1)) \cdot |S_1| \cdot I.
	\end{align*}
	Using Legendre's duplication formula for the Gamma function we have
	\begin{align}\label{eq I}
		I & = \frac{1}{(\sigma-1/2)^{2K-1}}\cdot \sqrt{\pi} \cdot \frac{\Gamma(K-1/2)}{\Gamma (K)} \notag\\
		& = \left(\frac{\L}{a}\right)^{2K-1} \cdot \pi \cdot 2^{2-2K}\cdot \binom{2K-2}{K-1}.
	\end{align}
	Furthermore, our assumption (ES $2K$) guarantees $\lim_{T \to \infty} |S_1|/N(T) \to 1$ as $a\to 0$.
	It follows that
	\begin{align*}
		\lim_{a \to 0^+} \lim_{T \to \infty}  \frac{(2a)^{2K-1}}{T (\log T)^{2K}}	\int_{\mathcal T_1} \left| \sum_{|s-\r|<\frac{c}{\L}} \frac{1}{s-\r} \right|^{2K} dt = \binom{2K-2}{K-1}.
	\end{align*}
	
	Next, we estimate $\int_{\mathcal T_2}$. It will be convenient (though not essential)  to describe $\mathcal T_2$ in a clearer way. To do so, let us think of a zero $\r$ colored green if $\r \in S_1$, and red if $\r \in [T, 2T] - S_1$. We separate the red zeros into clusters, where two red zeros belong to the same cluster if the distance between them is less than $2b/\L$. Therefore,  different clusters are at least $2b/\L$ apart. Call these clusters $C_1, ..., C_n$, and label zeros in $C_j$ as $\r_{j, 1}, ..., \r_{j, n_j}$ counting with multiplicity and ordered according to their heights. Thus, 
	\begin{align*}
		\int_{\mathcal T_2} \left| \sum_{|s-\r|<\frac{c}{\L}} \frac{1}{s-\r} \right|^{2K} dt  & = \sum_{j=1}^{n} \int_{\g_{j,1} - b/\L}^{\g_{j, n_j} + b/\L} \left| \sum_{|s-\r|<\frac{c}{\L}} \frac{1}{s-\r} \right|^{2K} dt.
	\end{align*}
	To treat the $2K$-th power of the absolute value, we use triangle inequality, expand the $2K$-th power as a sum, and then for each summand use the well-known inequality that for nonnegative $x_1, ..., x_n$ and positive integer $n$,
	\begin{align*}
		n\cdot	x_1 x_2\cdots x_n \le x_1^n + \cdots + x_n^n.
	\end{align*}
	In this way we see that
	\begin{align*}
		\left| \sum_{|s-\r|<\frac{c}{\L}} \frac{1}{s-\r} \right|^{2K}  & \le 	\left( \sum_{|s-\r|<\frac{c}{\L}} \frac{1}{|s-\r|} \right)^{2K} \\
		& \le \frac{1}{2K} \sum_{|s-\r|<\frac{c}{\L}} \frac{1}{|s-\r|^{2K}} \cdot (\text{ the number of times each summand appears}) \\
		& \le \frac{1}{2K} \sum_{|s-\r|<\frac{c}{\L}} \frac{1}{|s-\r|^{2K}} \cdot (2K \cdot \#(t)^{2K-1}),
	\end{align*}
	where $\#(t) = \sum_{|s-\r|<\frac{c}{\L}} 1$. It follows that
	\begin{align*}
		\int_{\mathcal T_2} \left| \sum_{|s-\r|<\frac{c}{\L}} \frac{1}{s-\r} \right|^{2K} dt  & \le \sum_{j=1}^{n} \int_{\g_{j,1} - b/\L}^{\g_{j, n_j} + b/\L} \sum_{|s-\r|<\frac{c}{\L}} \frac{1}{|s-\r|^{2K}} \cdot \#(t)^{2K-1} dt \\
		& \le \sum_{j=1}^{n}  \sum_{\r\in C_j} \int_{\g - b/\L}^{\g+ b/\L} \frac{1}{|s-\r|^{2K}} \cdot \#(t)^{2K-1} dt.
	\end{align*}
	Note that for $t \in (\g-b/\L, \g + b/\L)$ we have
	\begin{align*}
		\max_t \#(t) \le \text{ the number of zeros } \r^* \text{ with } |\r^*-\r|<2b/\L := \#_2(\r).
	\end{align*}
	Therefore,
	\begin{align}\label{eq int t2}
		\int_{\mathcal T_2} \left| \sum_{|s-\r|<\frac{c}{\L}} \frac{1}{s-\r} \right|^{2K} dt  & \le  \sum_{j=1}^{n}  \sum_{\r\in C_j} \#_2(\r)^{2K-1} \int_{\g - b/\L}^{\g+ b/\L} \frac{1}{|s-\r|^{2K}} dt \notag \\
		& \le   \sum_{j=1}^{n}  \sum_{\r\in C_j} \#_2(\r)^{2K-1} \int_{-\infty}^{\infty} \frac{1}{|s-\r|^{2K}} dt \notag \\
		& = I \cdot \sum_{\text{red } \r} \#_2(\r)^{2K-1}.
	\end{align}
	Now observe that
	\begin{align*}
		\sum_{\text{red } \r} \#_2(\r)^{2K-1} & = \sum_{\text{red } \r} (\text{ the number of } (2K-1) \text{-tuples in } (\g-2b/\L, \g+2b/\L)) \\
		& =  \sum_{\substack{\text{ all possible} (2K-1) \text{-tuples} \\ \text{appearing in the line above} } } \quad \sum_{\substack{\text{ red } \r \text{ with } |\g-\g^*|<2b/\L \\ \text{ for all } \g^* \text{ in the } (2K-1)\text{-tuple}}} 1 \\
		& = \text{ the number of } 2K \text{-tuples } (\lambda_1, \lambda_2 , ..., \lambda_{2K-1}, \r)
	\end{align*}
	where in the last line the $(\lambda_1, ..., \lambda_{2K-1})$ is a $(2K-1)$-tuple in the outer sum and $\r$ is a zero in the inner sum in the second last line. Separate the contribution of the $2K$-tuples in which all the $2K$ zeros are identically labelled in the list of zeta zeros, and we see that the above number is
	\begin{align*}
		\le \left (\sum_{\text{ red } \r} 1 \right ) + N(K, 4b, T),
	\end{align*}
	where we recall that $N(K,v, T)  $ is the number of $2K$-tuples $(\lambda_1, ..., \lambda_{2K}) $ of zeta zeros  in $[T, 2T]$, not all components having the same label in the list $\r_1, \r_2, ...$ of zeta zeros, and all components lying in an interval of length $v/\L$.
	
	Now the $2K$-tuple Essential Simplicity hypothesis guarantees that 
	\begin{align*}
		\lim_{T \to \infty} \frac{ \left (\sum_{\text{ red } \r} 1 \right ) + N(K, 4b, T)}{N(T)} \to 0 \quad \text{  as  } \; a\to 0. 
	\end{align*}
	This together with~\eqref{eq I} and~\eqref{eq int t2} gives
	\begin{align*}
		\lim_{a \to 0^+} \lim_{T \to \infty}  \frac{(2a)^{2K-1}}{T (\log T)^{2K}}		\int_{\mathcal T_2} \left| \sum_{|s-\r|<\frac{c}{\L}} \frac{1}{s-\r} \right|^{2K} dt  =0.
	\end{align*}
	\qed

	\section{Acknowledgments}
	
	I would like to thank the referee for a careful reading of the manuscript and for making a number of helpful suggestions.
	
	%%%%%%%%%%%%%%%%%%%%%%%%%%%%%%%%%%%%%%%%%%%%%%%%%

	%%%%%%%%%%%%%%%%%%%%%%%%%%%%%%%%%%%%%%%%%%%%%%%%%

\end{document}